\documentclass[12pt]{article}
\usepackage{latexsym,a4wide,amssymb,amsmath,graphicx,epsfig}

\newtheorem{thm}{Theorem}
\newtheorem{lem}{Lemma}
\newtheorem{cor}{Corollary}
\newtheorem{conj}{Conjecture}
\newtheorem{prop}{Proposition}
\newtheorem{exer}{Exercise}
\newtheorem{defi}{Definition}
\newtheorem{examp}{Example}

\newcommand{\ebox}{\hfill $\Box$\\\vspace{0.15cm}}
\newcommand{\pr}{{\bf Proof.}\ }
\newcommand{\bt}{\begin{thm}}
\newcommand{\et}{\end{thm}}
\newcommand{\bl}{\begin{lem}}
\newcommand{\el}{\end{lem}}
\newcommand{\bp}{\begin{prop}}
\newcommand{\ep}{\end{prop}}
\newcommand{\bc}{\begin{cor}}
\newcommand{\ec}{\end{cor}}
\newcommand{\bcj}{\begin{conj}}
\newcommand{\ecj}{\end{conj}}
\newcommand{\bex}{\begin{exer}}
\newcommand{\eex}{\end{exer}}
\newcommand{\bi}{\begin{itemize}}
\newcommand{\ei}{\end{itemize}}
\newcommand{\be}{\begin{equation}}
\newcommand{\ee}{\end{equation}}
\newcommand{\ben}{\begin{enumerate}}
\newcommand{\een}{\end{enumerate}}
\newcommand{\bd}{\begin{defi}}
\newcommand{\ed}{\end{defi}}
\newcommand{\bexa}{\begin{examp}}
\newcommand{\eexa}{\end{examp}}

\newcommand{\mt}{t\kern-0.035cm\char39\kern-0.03cm}
\newcommand{\ml}{l\kern-0.035cm\char39\kern-0.03cm}
\newcommand{\md}{d\kern-0.035cm\char39\kern-0.03cm}

\newcommand{\ovl}{\overline}
\newcommand{\noi}{\noindent}
\newcommand{\veps}{\varepsilon}

\begin{document}

\title{\vspace{-2.3cm} On the existence of aggregation functions with given super-additive and sub-additive
transformations}

\author{}
\date{}
\maketitle

\begin{center}
\vspace{-1.3cm}

{\large Alexandra \v{S}ipo\v{s}ov\'a, Ladislav \v Sipeky and Jozef \v{S}ir\'a\v{n} \footnote{Email addresses: $\{$siposova,sipeky,siran$\}$@math.sk}
\vspace{5mm}

{\small\it Faculty of Civil Engineering, Slovak University of Technology, 810 05 Bratislava, Slovakia}

}
\vspace{4mm}

\end{center}

\begin{abstract}
In this note we study restrictions on the recently introduced super-additive and sub-additive transformations, $A\mapsto A^*$ and $A\mapsto A_*$, of an aggregation function $A$. We prove that if $A^*$ has a slightly stronger property of being strictly directionally convex, then $A=A^*$ and $A_*$ is linear; dually, if $A_*$ is strictly directionally concave, then $A=A_*$ and $A^*$ is linear. This implies, for example, the existence of pairs of functions $f\le g$ sub-additive and super-additive on $[0,\infty[^n$, respectively, with zero value at the origin and satisfying relatively mild extra conditions, for which there exists no aggregation function $A$ on $[0,\infty[^n$ such that $A_*=f$ and $A^*=g$.

\vskip 2mm

\noi {\em Keywords: aggregation function, sub-additive and super-additive transformation}

\end{abstract}

\vskip 3mm

\section{Introduction}\label{intro}

A mapping $A: [ 0, \infty[^n \to [0, \infty[$ is called an aggregation function if $A({\bf 0})=A(0,\ldots,0)=0$ and $A$ is increasing in each coordinate. Literature on aggregation functions is abundant and we refer only to \cite{Bel,Grab} for basic facts. We will focus on a pair of mutually dual fundamental transformations of aggregation functions, which have been introduced in \cite{Greco}, motivated by interesting applications in economics.

\bd Let $A: [ 0, \infty[^n \to [0, \infty[$ be an aggregation function. The sub-additive transformation $A_*:[0,\infty [^n \to [0, \infty]$ of $A$ is given by

\begin{equation}
A_*({\bf x})= {\rm inf}\  \{\sum_{i=1}^k A({\bf x}^{(i)}) \ \  | \ \ \sum_{i=1}^k {\bf x}^{(i)} \geq {\bf x}\} \label{sub} \end{equation}

\noindent Similarly, the super-additive transformation $A^*:[0,\infty [^n \to [0, \infty]$ of $A$ is defined by

\begin{equation}
A^*({\bf x})= {\rm sup}\  \{\sum_{j=1}^{k} A({\bf x}^{(j)}) \ \ | \ \ \sum_{j=1}^{k} {\bf x}^{(j)} \leq {\bf x}\}\ . \label{super}\end{equation}
\ed

We note that the transformations (\ref{sub}) and (\ref{super}) were originally introduced in \cite{Greco} for aggregation functions as defined above {\sl and} with the property that $A_*$ and $A^*$ do not attain the value $\infty$.
\smallskip

It easy to show \cite{Greco} that the functions $A_*$ and $A^*$ are indeed, as their names suggest, sub-additive and super-additive, respectively, that is, $A_*({\bf u}+{\bf v})\le A_*({\bf u})+A_*({\bf v})$ and $A^*({\bf u}+{\bf v})\ge A^*({\bf u})+A^*({\bf v})$ for every ${\bf u},\ {\bf v}\in [0,\infty [^n$, where addition is defined coordinate-wise in the usual manner.
\smallskip

This suggests the question of whether or not for every pair $f,g:\ [0,\infty[^n\to [0,\infty]$ such that $f({\bf 0})=g({\bf 0})=0$, $f({\bf x})\le g({\bf x})$ for every ${\bf x}\in [0,\infty[^n$, with $f$ sub-additive and $g$ super-additive, there exists an aggregation function $A$ on $[0,\infty[^n$ such that $A_*=f$ and $A^*=g$.
\smallskip

In this paper we show that the answer to this question is negative if relatively mild extra conditions are imposed on $f$ and $g$. This is a consequence of our findings stemming from a more detailed study of replacing the super- and sub- additivity properties of the above transformations by the slightly stronger properties of directional convexity and concavity. Our main results say in a nutshell that if an aggregation function $A$ is such that $A^*$ is strictly directionally convex, then necessarily $A=A^*$ and $A_*$ is linear; dually, if $A_*$ is strictly directionally concave, then $A=A_*$ and $A^*$ is linear.
\smallskip

We explain our tools in section \ref{1dim} on the one-dimensional case first, as it clearly distinguishes what hurdles one needs to overcome in extending the result to the multi-dimensional case, which is done in section \ref{ndim}. We conclude in section \ref{con} by a discussion.

\section{The one-dimensional case}\label{1dim}

We consider here the one-dimensional case, a preliminary report on which can be found in \cite{SiSi}. We recall that a function $h:\ [0,\infty[\to [0,\infty[$ is strictly convex if for every distinct $u,v\in [0,\infty[$ and every $t\in\ ]0,1[$ we have $h(tu+(1-t)v) < th(u) + (1-t)h(v)$; the usual concept of convexity s obtained by replacing the strict inequality with a non-strict one. A function $h:\ [0,\infty[\to [0,\infty[$ is called {\sl directionally convex} if its increments are non-decreasing, that is, if $h(x+t)-h(x)\le h(y+t)-h(y)$ whenever $0\le x\le y$ and $h\ge 0$. It is easy to see that this condition is equivalent to $h(x)+h(y) \le h(u)+h(v)$ for every $u,v,x,y\in [0,\infty[$ with $u\le x,y\le v$ and $u+v=x+y$; we will call $h$ {\sl strictly directionally convex} if all the inequalities in the are replaced with strict ones. Recall also that $h$ is {\sl strictly super-additive} if $h(x)+h(y) < h(x+y)$ for every distinct $x,y\in [0,\infty[$; allowing non-strict inequality gives back the concept of super-additivity.
\smallskip

As we shall see in the next lemma, directional convexity is equivalent to convexity for functions locally bounded at some point on the real axis. We prefer, however, using directional convexity in the statement of our main result of this section, because this turns out to be the way to extend our findings to the multi-dimensional case later. We will restrict our attention here to functions defined on $[0,\infty[$, as only these are of our concern, and we will say that such a function $h$ is {\sl locally bounded at some point} if there exists an $x_0\in]0,\infty[$ such that $h$ is bounded in some open neighbourhood of $x_0$.

\bl\label{lem1}
Let $h:\ [0,\infty[\to [0,\infty[$ be a function locally bounded at some point. Then:
\smallskip

{\rm (a)} The function $h$ is directionally convex if and only if it is convex, and $h$ is strictly directionally convex if and only if it is strictly convex; moreover, any of these properties imply continuity of $h$ on $[0,\infty[$.
\smallskip

{\rm (b)} If $h$ is directionally convex and $h(0)=0$, then $h$ is super-additive, and if $h$ is strictly directionally convex with $h(0)=0$, then it is strictly super-additive.
\smallskip

{\rm (c)} If $h$ is directionally convex and $h(0)=0$, then for every $u,v\in ]0,\infty[$ such that $0<u<v$ we have $h(u)/u\le h(v)/v$; in particular, $h(\veps)\le \veps h(1)$ for every $\veps\in ]0,1]$.
\el

\pr
For (a), taking $x=y=(u+v)/2$ in the definition of directional convexity of $h$ implies that $h(\frac{1}{2}(u+v)) \le \frac{1}{2}h(u)+\frac{1}{2}h(v)$; this property is known as midpoint convexity. Since $h$ is assumed to be locally bounded at some point, midpoint convexity of $h$ implies its convexity and continuity by the Bernstein-Doetsch theorem \cite{Ber}. This shows that under our assumptions, directional convexity of $h$ implies its convexity. Conversely, if $u\le x,y\le v$ are such that $u+v=x+y$, then there exists a $t\in\ [0,1]$ such that $x=tu+(1-t)v$ and $y=(1-t)u+tv$. Applying convexity of $h$ gives $h(x)\le th(u)+(1-t)h(v)$ and $h(y)\le (1-t)h(u)+th(v)$, and adding the last two inequalities together gives $h(x)+h(y)\le h(u)+h(v)$. It follows that convexity of $h$ implies its directional convexity. The ``strict'' versions of the two implications are obvious.
\smallskip

Item (b) follows from directional convexity of $h$ with $h(0)=0$ by letting $u=0$ and $v=x+y$ and using strict inequalities if $h$ is strictly directionally convex.
\smallskip

To prove (c) it is sufficient to apply (a), and hence convexity of $h$, to the equation $u=(1-t)0+tv$ for $t=u/v$, by which $h(u)\le (1-t)h(0)+th(v)=uh(v)/v$; the last part of (c) follows by letting $u=\veps$ and $v=1$.
\ebox

We note that if the assumption of $h$ being locally bounded at some point on is dropped, then directional convexity need not imply convexity. Counterexamples constructed with the help of Hamel bases of the vector space of the real numbers over the field of rational numbers are identical to those showing that midpoint convexity need not imply convexity; we refer to \cite{Kuc} for details together with accompanying theory.
\smallskip

We are now ready to prove our main result in the one-dimensional case and we will do so in the directional convexity setting.

\bt\label{main1}
Let $A:\ [0,\infty[ \to [0,\infty]$ be an aggregation function. If $A^*$ is strictly directionally convex, then $A(x)=A^*(x)$ for every $x\in [0,\infty[$. Moreover, in such a case we have $A_*(x)=cx$ for every $x\in [0,\infty[$, where $c=\lim_{t\to 0^+} A(t)/t$.
\et

\pr
Let $A:\ [0,\infty[ \to [0,\infty]$ be an aggregation function such that $A^*$ is strictly directionally convex; note that this automatically implies that $A^*$ is locally bounded at some (in fact, at every) point in $]0,\infty[$. Clearly,  $A(x)\le A^*(x)$ and $A^*(x)>0$ for every $x\in\ ]0,\infty[$, and $A(0)=A^*(0)=0$; moreover, both $A$ and $A^*$ are increasing functions. Suppose for a contradiction that there exists an $\bar x>0$ such that $A(\bar x) <A^*(\bar x)$, where, using an equivalent form of the definition of a super-additive transformation,
\begin{equation}\label{nula}
A^*(\bar x)={\rm sup}\left\{\sum_{j=1}^{k} A(x_j)\ |\ x_1{\ge}x_2{\ge}{\ldots}{\ge}x_{k}{>}0\ \ {\rm and}\ \ \sum_{j=1}^{k} x_j=\bar x\right\} \ .
\end{equation}
If the sums in (\ref{nula}) consist of just one element, $A(\bar x)$ and $\bar x$, then we have $A(\bar x) < A^*(\bar x)$ by our assumption, and so we may let $k\ge 2$ throughout. Let $\delta_1=(A^*(\bar x)-A(\bar x))/2 > 0$ and define $\veps>0$ by $\veps = {\rm min}\ \{\ \delta_1/A^*(1)\ ,\ \bar x/3\ ,\ 1\ \}$. We will distinguish two cases.
\smallskip

{\sl Case 1:} \ $x_1\in\ ]\bar x-\veps,\bar x[$, so that $\sum_{j=2}^k < \veps$. Since $A$ is increasing, we have $A(x_1) < A(\bar x)$, and the way $\delta_1$ was introduced then implies that $A(x_1) < A(\bar x)= A^*(\bar x)-2\delta_1$. Further, by earlier inequalities, the facts that $A^*$ is super-additive and increasing, and by the last statement of part (c) of Lemma \ref{lem1} we have
\[\sum_{j=2}^k A(x_j) \le \sum_{j=2}^k A^*(x_j) \le A^*(\sum_{j=2}^k x_j)\le A^*(\veps)\le \veps A^*(1)\ .\]
These findings together with $\veps \le \delta_1/A^*(1)$ lead to the estimate
\be\label{est1}
\sum_{j=1}^k A(x_j)=A(x_1) + \sum_{j=2}^k A(x_j)\le A^*(\bar x)-2\delta_1 + \veps A^*(1) \le A^*(\bar x) - \delta_1\ .
\ee

{\sl Case 2:} $x_1\le \bar x-\veps$. Let $\delta_2=A^*(\bar x)-A^*(\bar x-\veps)-A^*(\veps)$; by {\sl strict} super-additivity of $A^*$ (a consequence of strict directional convexity by parts (a) and (b) of Lemma \ref{lem1}) we have $\delta_2>0$. Let $\ell$ be the smallest subscript for which $\sum_{j=1}^\ell x_j \in [\veps, \bar x-\veps]$; note that $\ell$ is well defined and $1\le \ell\le k-1$, since $x_j\le x_1$ ($1 \le j \le k$) and $\veps < \bar x/3$ by our choice of $m$. Let $u=\sum_{j=1}^\ell x_j$ and $v=\sum_{j=\ell+1}^k x_j$, with $u+v=\bar x$. The facts that $A^*$ is increasing and super-additive now imply
\[  \sum_{j=1}^k A(x_j)=\sum_{j=1}^\ell A(x_j) + \sum_{j=\ell+1}^k A(x_j)\le \sum_{j=1}^\ell A^*(x_j) + \sum_{j=\ell+1}^k A^*(x_j) \le  A^*(u) + A^*(v)\ .\]
We proceed by applying (not necessarily strict) directional convexity of $A^*$ to the four values $u\ge \veps$ and $v\le\bar x-\veps$,  by which $A^*(u)+A^*(v)\le A^*(\bar x-\veps)+A^*(\veps)$. Combining this with the previous inequality in conjunction with the definition of $\delta_2$ we obtain
\be\label{est2}
\sum_{j=1}^k A(x_j) \le A^*(\bar x-\veps)+A^*(\veps) = A^*(\bar x) - \delta_2\ .
\ee

But now, letting $\delta={\rm min}\{\delta_1,\delta_2\} >0$ one sees that by (\ref{est1}) and (\ref{est2}) we have in (\ref{nula}) the inequality $\sum_{j=1}^{k}A(x_j) \le A^*(\bar x) - \delta$ whenever $k\ge 2$, and we know that $A(\bar x) < A^*(\bar x)$. Therefore, by (\ref{nula}) we arrive at the absurd conclusion that $A^*(\bar x) < A^*(\bar x)$. This contradiction completes the proof of the fact that $A(x)=A^*(x)$ for every $x\in [0,\infty[$.
\smallskip

It follows that $A$ (being equal to $A^*$) is convex. By part (b) of Lemma \ref{lem1} the function $A(t)/t$ is increasing on $]0,\infty[$ and so the limit of $A(t)/t$ as $t\to 0^+$ exists. If $c=\lim_{t\to 0^+} A(t)/t$, we have $A_*(x)=cx$ by Corollary 1 of \cite{Sipo}.
\ebox

Note that no part of the proof refers to any other assumptions that one often makes about functions, in particular, to continuity, although our assumptions on $A^*$ imply that $A^*$ is continuous, by the Bernstein-Doetsch theorem \cite{Ber}; cf. part (a) of Lemma \ref{lem1}.
\smallskip

By the obvious duality we have the following result the proof of which is left to the reader. To have a corresponding companion to the notion of strict directional convexity, a function $h:\ [0,\infty[\to [0,\infty[$ is said to be {\sl strictly directionally concave} if for every distinct $u,v,x,y\in [0,\infty[$ such that $u < x,y<v$ and $u+v=x+y$ we have $h(u) + h(v) < h(x) + h(y)$.

\bt\label{main2} Let $A:\ [0,\infty[ \to [0,\infty]$ be an aggregation function. If $A_*$ is strictly directionally concave, then $A(x)=A_*(x)$ for every $x\in [0,\infty[$. Moreover, in such a case we have $A^*(x)=cx$ for every $x\in [0,\infty[$, where $c=\lim_{t\to 0^+} A(t)/t$.  \hfill $\Box$
\et

A contrapositive version of the two theorems gives, in general, a negative answer to the question of existence of one-dimensional aggregation functions with arbitrary preassigned sub-additive and super-additive transformations.

\bt\label{main3}
Let $f,g:\ [0,\infty[ \to [0,\infty]$ be such that $f(0)=g(0)=0$ and $f(x)\le g(x)$ for every $x\in [0,\infty[$. If $f$ is strictly directionally concave and $g$ is super-additive but not linear, or else if $f$ is sub-additive but not linear and $g$ is strictly directionally convex, then there is no aggregation function $A$ on $[0,\infty[$ such that $A_*(x)=f(x)$ and $A^*(x)=g(x)$ for every $x\in [0,\infty[$. \hfill $\Box$
\et

\section{The multi-dimensional case}\label{ndim}

We will use the standard notation ${\bf x}=(x_1,x_2,\ldots,x_n)$, ${\bf y}=(y_1,y_2,\ldots,y_n)$, and so on, for points in $[0,\infty[^n$; in particular, ${\bf 0}$ and ${\bf 1}$ stand for the points $(0,0,\ldots,0)$ and $(1,1,\ldots,1)$. We will write ${\bf x}\le {\bf y}$ if ${\bf y}-{\bf x} \in [0,\infty[^n$, and ${\bf x}< {\bf y}$ if ${\bf x}\le {\bf y}$ but ${\bf x}\ne {\bf y}$. As usual, ${\bf e}_i$ will denote the $i$-th unit vector, so that for every ${\bf x}=(x_1,x_2,\ldots,x_n)\in [0,\infty[^n$ we have ${\bf x}=\sum_{i=1}^nx_i{\bf e}_i$.
\smallskip

Extending the concept of directional convexity introduced in section \ref{1dim} to the multidimensional case, we will say that a function $h:\ [0,\infty[^n\to [0,\infty[$ is {\sl directionally convex} if $h({\bf x})+h({\bf y}) \le h({\bf u})+h({\bf v})$ whenever ${\bf u},{\bf v},{\bf x},{\bf y}\in [0,\infty[^n$ are such that ${\bf u}\le {\bf x},{\bf y}\le {\bf v}$ and ${\bf u}+{\bf v}={\bf x}+{\bf y}$. Moreover, such a function $h$ will be said to be {\sl strictly directionally convex} if all three inequalities in the preceding definition are strict. For an equivalent definition with the help of increments (mentioned in the one-dimensional case) we refer, for instance, to \cite{ArSc}. An alternative term for directional convexity is ultra-modularity, cf. \cite{MaMo}. We also note without going into details that directional convexity of an aggregation function is equivalent to its super-modularity (see e.g. \cite{ArSc}) and coordinate-wise convexity.
\smallskip

The relationship between directional convexity, super-additivity and convexity is more delicate in two or more dimensions compared with the one-dimensional case. Clearly, if $h$ is a function as above such that $h({\bf 0})=0$, then directional convexity of $h$ implies its super-additivity. Directional convexity of $h$ also implies its midpoint convexity in each coordinate, meaning that for every fixed $i\in \{1,2,\ldots, n\}$ and every $(n-1)$-tuple of non-negative real numbers $x_1,\ldots,x_{i-1},x_{i+1},\ldots,x_n$ the function $z\mapsto h(x_1,\ldots,x_{i-1}, z , x_{i+1},\ldots,x_n)$ for $z\in [0,\infty[$ is midpoint convex. To see this it is sufficient to consider four points ${\bf u}, {\bf v}, {\bf x},{\bf y}$ as above such that $u_j=x_j$ and $v_j=y_j$ for and all $j\ne $, $1\le j\le n$ and $x_i=y_i=(u_i+v_i)/2$.
\smallskip

Hence, by Lemma \ref{lem1}, if every restriction of $h$ to one coordinate (fixing the value of other $n-1$ coordinates arbitrarily) is locally bounded at some point, then $h$ is continuous when restricted to an arbitrary coordinate, and also coordinate-wise convex. Recall that he latter means that or every $i\in\{1,2,\ldots,n\}$ and for every pair of ${\bf x},{\bf y}\in [0,\infty[^n$ such that ${\bf x}-{\bf y}$ is a scalar multiple of ${\bf e}_i$ we have $h(t{\bf x}+(1-t){\bf y}) \le th({\bf x}) + (1-t)h({\bf y})$ for every $t\in [0,1]$. However, directional convexity in two or more dimensions does {\sl not} imply convexity, as one may see for $n=2$ on a trivial example of an aggregation function $h:\ [0,\infty[^2\to [0,\infty[$ given by $h(x_1,x_2)=(x_1+1)(x_2+1)-1$. For completeness, note also that convexity of an aggregation function in two or more dimensions does not imply super-additivity; consider, for instance, $n=2$ and $h:\ [0,\infty[^2\to [0,\infty[$ given by $h(x_1,x_2)= (x_1-x_2)^2 +4x_2^2$.
\smallskip

We will summarize and slightly extend the above discussion in the form of an auxiliary result. To avoid using multidimensional versions of the Bernstein-Doetsch theorem (see e.g. \cite{Mur}) we replace the local boundedness  assumption mentioned above with a global one as we will only be dealing with `tame' functions later in our main results.

\bl\label{lem2}
Let $h:\ [0,\infty[^n\to [0,\infty[$ be a function bounded on every bounded subset of $[0,\infty[^n$. Then:
\smallskip

{\rm (a)} Directional convexity of $h$ implies its coordinate-wise convexity and continuity in every coordinate.
\smallskip

{\rm (b)} If $h$ is directionally convex and $h({\bf 0})=0$, then $h$ is super-additive, and if $h$ is strictly directionally convex with $h({\bf 0})=0$, then it is strictly super-additive.
\smallskip

{\rm (c)} If $h$ is directionally convex and $h({\bf 0})=0$, then $h(\veps{\bf 1})\le \veps h({\bf 1})$ for every $\veps\in [0,1]$.
\el

\pr
Items (a) and (b) have been proved in the discussion preceding the lemma. Regarding (c), let $f(x)=h(x,x,\ldots, x)$ for every $x\in [0,\infty[$. Since $f(0)=0$ and $f$ is obviously directionally convex as a function of one variable, part (c) of Lemma \ref{lem1} implies that $f(\veps)\le \veps f(1)$, which means that $h(\veps{\bf 1})\le \veps h({\bf 1})$ for every $\veps\in [0,1]$. \ebox

We are now ready to state and prove the main result of this section.

\bt\label{main-1}
Let $A:\ [0,\infty[^n \to [0,\infty]$ be an aggregation function. If $A^*$ is strictly directionally convex, then $A({\bf x})=A^*({\bf x})$ for every ${\bf x}\in [0,\infty[^n$. Moreover, in such a case we have $A_*({\bf x})={\bf \nabla \cdot x}$ for every ${\bf x}\in [0,\infty[^n$, where $\nabla_i=\lim_{t\to 0^+} A(t{\bf e_i})/t$.
\et

\pr
Let $A:\ [0,\infty[^n \to [0,\infty]$ be an aggregation function such that $A^*$ is strictly directionally convex. We obviously have $A({\bf x})\le A^*({\bf x})<\infty$, the second inequality being a consequence of strict directional convexity. We also have $A^*({\bf x})>0$ for every ${\bf x}\in\ [0,\infty[^n{\setminus}\{\bf 0\}$, and $A({\bf 0})=A^*({\bf 0})=0$; moreover, both $A$ and $A^*$ are increasing in every coordinate. Note that super-additivity (a consequence of directional convexity, by Lemma \ref{lem2}) of $A^*$ implies its boundedness on every bounded set.
\smallskip

Suppose for a contradiction that there exists an $\ovl{\bf x}\ne {\bf 0}$ such that $A(\bar x) <A^*(\bar x)$, where, using an equivalent form of the definition of a super-additive transformation,
\begin{equation}\label{e-1}
A^*(\ovl{\bf x})={\rm sup}\left\{\ \sum_{j=1}^{k} A({\bf x}^{(j)})\ |\ {\bf 0}\ne {\bf x}^{(j)} \in [0,\infty[^n \  (1\le j\le k), \  \sum_{j=1}^{k} {\bf x}^{(j)}=\ovl{\bf x}\ \right\} \ .
\end{equation}
If both sums consist of just one element, $A(\ovl{\bf x})$ and $\ovl{\bf x}$, then we have $A(\ovl{\bf x}) < A^*(\ovl{\bf x})$, and so we will assume that $k\ge 2$ in forthcoming considerations. Also, we may assume that all coordinates of ${\ovl{\bf x}}$ are positive. Indeed, if, say, without loss of generality, ${\bar x}_n=0$, then we would have $x^{(j)}_n=0$ for all $j\in \{1,2,\ldots,k\}$ in the points entering (\ref{e-1}), which just means reducing the dimension from $n$ to $n-1$.
\smallskip

Before we proceed we need to introduce several parameters. Let $\xi>0$ be the smallest value of the coordinates $\ovl{x}_i$ of the point $\ovl{\bf x}$. Further, let $\delta_1>0$ be defined by $2\delta_1=A^*(\ovl{\bf x}) - A(\ovl{\bf x})$. For every {\sl proper} subset $I\subset \{1,2,\ldots,n\}$ and every $\eta\ge 0$ we let ${\ovl{\bf x}}_{I,\eta}$ denote the point whose $i$-th coordinate is equal to $\bar x_i$ for every $i\in I$ and to $\eta$ for every $i\in J=\{1,\ldots,n\}{\setminus} I$. Now, for every such {\sl non-trivial} partition $\{I,J\}$ of $\{1,2,\ldots,n\}$, strict directional convexity of $A^*$ implies that $A^*(\ovl{\bf x}_{I,0}) + A^*(\ovl{\bf x}_{J,0}) < A^*(\ovl{\bf x})$. This strict inequality together with the fact that $A^*$ is coordinate-wise continuous by Lemma \ref{lem2} imply that there exist $\mu >0$ and $\delta_2>0$ such that for every $\nu\in [0,\mu]$ and every non-trivial partition $\{I,J\}$ of $\{1,2,\ldots,n\}$ we have
\begin{equation}\label{e-1.5}
A^*(\ovl{\bf x}_{I,\nu}) + A^*(\ovl{\bf x}_{J,\nu}) < A^*(\ovl{\bf x})-\delta_2\ .
\end{equation}
Finally, we introduce $\veps>0$ by letting
\begin{equation}\label{e-2}
\veps = {\rm min}\ \left\{ \  \frac{\xi}{3}\ ,\ \frac{\delta_1}{A^*({\bf 1})}\ ,\
\mu\ ,\ 1  \ \right\}\ .
\end{equation}
We will distinguish three cases.
\smallskip

{\sl Case 1:} Suppose that in the $k$-tuple ${\bf x}^{(j)}$ ($1\le j\le k$) entering (\ref{e-1}) there exists a superscript $j\in\{1,2,\ldots,k\}$ such that ${\bf x}^{(j)}\ge {\ovl{\bf x}}-\veps{\bf 1}$. Without loss of generality we may assume that $j=1$ and then, obviously, $\sum_{j=2}^{k}{\bf x}^{(j)} \le \veps{\bf 1}$.
\smallskip

In order to estimate $\sum_{j=1}^k A({\bf x}^{(j)})$ in the Case 1 we make a few observations. By our choice of $\delta_1$ and the fact that $A$ is increasing in every coordinate we have $A({\bf x}^{(1)})\le A(\ovl{\bf x}) = A^*(\ovl{\bf x}) -2\delta_1$. Moreover, the facts that $A^*$ is increasing and super-additive applied to the inequality $\sum_{j=2}^{k}{\bf x}^{(j)} \le \veps{\bf 1}$ implies that $\sum_{j=2}^k A({\bf x}^{(j)}) \le \sum_{j=2}^k A^*({\bf x}^{(j)}) \le A^*(\veps{\bf 1})$. Applying part (c) of Lemma \ref{lem2} we obtain $A^*(\veps{\bf 1})\le \veps A^*({\bf 1})$. Putting the pieces together and using our choice of $\veps\le \delta_1/A^*({\bf 1})$ we arrive at our first partial conclusion: If $\sum_{j=1}^{k} {\bf x}^{(j)}$ is as in the Case 1, then
\begin{equation}\label{e-3}
\sum_{j=1}^k A({\bf x}^{(j)}) = A({\bf x}^{(1)}) + \sum_{j=2}^k A({\bf x}^{(j)}) \le A^*(\ovl{\bf x}) -2\delta_1 + \veps A^*({\bf 1}) \le A^*(\ovl{\bf x}) -\delta_1\ .
\end{equation}

{\sl Case 2:} Suppose that the $k$-tuple ${\bf x}^{(j)}$ ($1\le j\le k$) appearing in (\ref{e-1}) has the property that for every $i\in \{1,2,\ldots,n\}$ there exists a $j=j_i\in \{1,2,\ldots, k\}$ such that $x^{(j)}_i \ge {\bar x}_i- \veps$ but we do {\sl not} have ${\bf x}^{(j)}\ge {\ovl{\bf x}}-\veps{\bf 1}$. We thus may suppose without loss of generality that there is an $r\in \{1,\ldots,n-1\}$ for which $j_1=\ldots=j_r=1$ but $j_i\ge 2$ for all $i$ such that $r+1\le i\le n$.
\smallskip

Let ${\bf y}={\bf x}^{(1)}$ and ${\bf z}=\sum_{j=2}^k{\bf x}^{(j)}$. With the help of the dominance of $A$ by $A^*$ and super-additivity of $A^*$ we obtain
\begin{equation}\label{e-4}
\sum_{j=1}^k A({\bf x}^{(j)}) \le  A^*({\bf x}^{(1)}) + \sum_{j=2}^k A^*({\bf x}^{(j)}) \le A^*({\bf x}^{(1)}) + A^*( \sum_{j=2}^k {\bf x}^{(j)}) = A^*({\bf y}) + A^*({\bf z}) \ .
\end{equation}
For the partition $\{I,J\}$ of $\{1,2,\ldots,n\}$ given by $I=\{1,\ldots,r\}$ and $J=\{r+1,\ldots,n\}$, let $\ovl{\bf x}_{I,\veps}$ and $\ovl{\bf x}_{J,\veps}$ be points as introduced earlier when defining the values of $\delta_2$ and $\veps$. By the assumptions of Case 2 and the way $r$ has been introduced we have ${\bf y}\le \ovl{\bf x}_{I,\veps}$ and ${\bf z}\le \ovl{\bf x}_{J,\veps}$. But by (\ref{e-1.5}) we know that $A^*(\ovl{\bf x}_{I,\veps}) + A^*(\ovl{\bf x}_{J,\veps}) \le A^*(\ovl{\bf x})-\delta_2$. Coupling these inequalities with (\ref{e-4}) yields our second partial conclusion: If $\sum_{j=1}^{k} {\bf x}^{(j)}$ is as in the Case 2, then
\begin{equation}\label{e-5}
\sum_{j=1}^k A({\bf x}^{(j)}) \le A^*({\bf y}) + A^*({\bf z}) \le  A^*({\ovl{\bf x}}_{I,\veps}) + A^*({\ovl{\bf x}}_{J,\veps}) \le  A^*({\ovl{\bf x}}) - \delta_2\ .
\end{equation}

{\sl Case 3:} Suppose that the $k$-tuple ${\bf x}^{(j)}$ ($1\le j\le k$) entering (\ref{e-1}) has the property that there exists an $i\in \{1,\ldots,n\}$ such that for each $j\in \{1,\ldots,k\}$ we have $x^{(j)}_i \le {\bar x}_i- \veps$. We may, of course, assume that $i=1$ and that the notation is chosen in such a way that $x^{(1)}_1\ge x^{(2)}_1\ge \ldots \ge x^{(k)}_1$. Let $\ell$ be the smallest index such that $\sum_{j=1}^\ell x^{(j)}_1 \ge \veps$. Our Case 3 assumptions imply that $1\le \ell \le n-1$ and that both sums $\sum_{j=1}^\ell x^{(j)}_1$ and $\sum_{j=\ell+1}^n x^{(j)}_1$ are in the interval $[\veps,\bar x_1-\veps]$.
\smallskip

Let ${\bf x}=\sum_{j=1}^\ell {\bf x}^{(j)}$ and ${\bf y}=\sum_{j=\ell+1}^n {\bf x}^{(j)}$, and let ${\bf u}=\veps{\bf e}_1$ and ${\bf v}={\ovl{\bf x}}-\veps{\bf e}_1$. It is clear that ${\bf x}+{\bf y} = {\bf u} + {\bf v} = {\ovl{\bf x}}$, and by our choice of $\ell$ we have ${\bf u} \le {\bf x}, {\bf y} \le {\bf v}$. Applying the assumed directional convexity to the four points just introduced gives $A^*({\bf x}) + A^*({\bf y}) \le A^*({\bf u}) + A^*({\bf v})$, and the sum on the right-hand side is always away from $A^*({\ovl{\bf x}})$. To formalize this, let $\delta_3 = {\rm min}\ \{ A^*(\ovl{\bf x}) - A^*(\veps{\bf e}_i) - A^*({\ovl{\bf x}}-\veps{\bf e}_i);\ 1\le i\le n\}$; observe that $\delta_3>0$ because of our choice of $\veps$ and our assumption that $\bar x_i>0$, $1\le i\le n$. The inequalities we have obtained in passing imply that if $\sum_{j=1}^{k} {\bf x}^{(j)}$ is as in the Case 3, then, for some $i$ (which was assumed to be $1$ above),
\begin{equation}\label{e-6}
\sum_{j=1}^k A({\bf x}^{(j)}) \le A^*({\bf x}) + A^*({\bf y}) \le A^*(\veps{\bf e}_i) + A^*({\ovl{\bf x}}-\veps{\bf e}_i) \le A^*(\ovl{\bf x}) - \delta_3\ .
\end{equation}

The conclusion regarding $A^*$ is now immediate. It is clear that for every $k\ge 2$ a $k$-tuple ${\bf x}^{(j)}$, $1\le j\le k$, such that $\sum_{j=1}^{k} {\bf x}^{(j)}=\ovl{\bf x}$, falls under one of the three cases considered above. Moreover, letting $\delta={\rm min}\{\delta_1,\delta_2,\delta_3\} >0$ one sees that by (\ref{e-3}), (\ref{e-5}) and (\ref{e-6}) we have in (\ref{e-1}) the inequality $\sum_{j=1}^{k}A({\bf x}^{(j)}) \le A^*(\ovl{\bf x}) - \delta$ whenever $k\ge 2$, and we know that $A(\ovl{\bf x}) < A^*(\ovl{\bf x})$. By (\ref{e-1}) we thus have $A^*(\ovl{\bf x}) < A^*(\ovl{\bf x})$, a contradiction. This proves that $A({\bf x})=A^*({\bf x})$ for every ${\bf x}\in [0,\infty[^n$.
\smallskip

It remains to deal with the statement on $A_*$. By part (a) of Lemma \ref{lem2} applied to $A=A^*$ one sees that $A$ is coordinate-wise convex. Consider an arbitrary $i\in\{1,2,\ldots,n\}$. By part (b) of Lemma \ref{lem1} the function $t\mapsto A(t{\bf e}_i)/t$ is increasing on $]0,\infty[$ and so the limit of $A(t{\bf e}_i)/t$ as $t\to 0^+$ exists. Let ${\bf \nabla}$ be the vector with $\nabla_i=\lim_{t\to 0^+} A(t{\bf e}_i)/t$. Applying Theorem 1 of \cite{Sipo} to the function $x_i\mapsto A(x_i{\bf e}_i)$ of one variable $x_i\in [0,\infty[$ we obtain the inequality $A_*(x_i{\bf e}_i) \le \nabla_ix_i$ for every $x_i\in [0,\infty[$. Since $A_*$ is sub-additive \cite{Greco}, for every ${\bf x}\in [0,\infty[^n$ we have
\begin{equation}\label{e-7} A_*({\bf x})=A_*(\sum_{i=1}^n x_i{\bf e}_i) \le \sum_{i=1}^n A_*(x_i{\bf e}_i) \le {\bf \nabla\cdot x}\ .
\end{equation}

We now prove the reverse inequality. By super-additivity of $A^*=A$ and the inequality $A^*(x_i{\bf e}_i) \ge \delta_ix_i$, which again is a consequence of Theorem 1 of \cite{Sipo} applied to the function $x_i\mapsto A(x_i{\bf e}_i)$, we obtain
\begin{equation}\label{e-8}
A({\bf x}) =A^*({\bf x}) = A^*(\sum_{i=1}^k x_i{\bf e}_i) \ge \sum_{i=1}^n A^*(x_i{\bf e}_i) \ge \sum_{i=1}^k \nabla_ix_i = {\bf \nabla\cdot x}
\end{equation}
for every ${\bf x}\in [0,\infty[^n$. But (\ref{e-8}) implies that $A_*({\bf x}) \ge ({\bf \nabla\cdot x})_*$ and since for linear functions of ${\bf x}\in [0,\infty[^n$ such as ${\bf \nabla\cdot x}$ we have $({\bf \nabla\cdot x})_*={\bf \nabla\cdot x}$, it follows that $A_*({\bf x}) \ge {\bf \nabla\cdot x}$ for every ${\bf x}\in [0,\infty[^n$. This together with (\ref{e-7}) implies that $A_*({\bf x})={\bf \nabla \cdot x}$ for every ${\bf x}\in [0,\infty[^n$, completing the proof. \ebox

The reader may have noticed that it is the case 2 of the above proof which covers a situation that does not appear in the one-dimensional case, while handling the cases 1 and 3 is an extension of the way the corresponding instances have been treated previously.
\smallskip

A dual version of Theorem \ref{main-1} is obtained in an obvious way. A function $h:\ [0,\infty[^n\to [0,\infty[$ is said to be {\sl strictly directionally concave} if for every distinct ${\bf u},{\bf v},{\bf x},{\bf y}\in [0,\infty[^n$ such that ${\bf u} < {\bf x},{\bf y}<{\bf v}$ and ${\bf u}+{\bf v}={\bf x}+{\bf y}$ we have $h({\bf u}) + h({\bf v}) < h({\bf x}) + h({\bf y})$. An appropriate modification of the above proof gives:

\bt\label{main-2} Let $A:\ [0,\infty[^n \to [0,\infty]$ be an aggregation function. If $A^*$ is strictly directionally concave, then $A({\bf x})=A_*({\bf x})$ for every ${\bf x}\in [0,\infty[^n$. Moreover, in such a case we have $A^*({\bf x})={\bf \nabla \cdot x}$ for every ${\bf x}\in [0,\infty[^n$, where $\nabla_i=\lim_{t\to 0^+} A(t{\bf e_i})/t$.  \hfill $\Box$
\et

Again, we state a contrapositive version of the two theorems which, in general, provide a negative answer to the question of existence of multi-dimensional aggregation functions with arbitrary preassigned sub-additive and super-additive transformations.

\bt\label{main-3}
Let $f,g:\ [0,\infty[^n \to [0,\infty]$ be such that $f({\bf 0})=g({\bf 0})=0$ and $f({\bf x})\le g({\bf x})$ for every ${\bf x}\in [0,\infty[^n$. If $f$ is strictly directionally concave and $g$ is super-additive but not linear, or else if $f$ is sub-additive but not linear and $g$ is strictly directionally convex, then there is no aggregation function $A$ on $[0,\infty[^n$ such that $A_*({\bf x})=f({\bf x})$ and $A^*({\bf x})=g({\bf x})$ for every ${\bf x}\in [0,\infty[^n$. \hfill $\Box$
\et

\section{Concluding remarks}\label{con}

The question of whether (or what kind of) further relaxations of assumptions in our main results can be made remains open. Here we just point out that, in general, the assumptions of {\sl strict} directional convexity or concavity cannot be completely dropped. We illustrate this by an example in dimension 1 where, as we saw in section \ref{1dim}, strict directional convexity is equivalent to strict convexity.
\smallskip

\bexa\label{ex1} Consider a function $A:\ [0,\infty[ \to [0,\infty[$ defined by
\begin{center}
$A(x) = \begin{cases}
x, & \mbox{if } x\in [0,4]; \\
\frac{1}{2}x+2, & \mbox{if } x\in [4,6]; \\
\frac{5}{6}x, & \mbox{if } x\in [6,12]; \\
\frac{5}{4}x-5 & \mbox{if } x\in [12,\infty[. \end{cases}$
\end{center}
\noindent Further, let us introduce functions $f,g:\ [0\infty[\to [0,\infty[$ given by
\begin{center}
$\begin{array}{rl}
f(x) = \begin{cases} x & \mbox{if } x\in [0,4]; \\ \frac{1}{2}x+2 & \mbox{if } x\in [4,6]; \\
\frac{5}{6}x & \mbox{if } x\in [6,\infty[; \end{cases}  \ \ \ &
\ \ \  g(x) = \begin{cases} x & \mbox{if } x\in [0,20]; \\ \frac{5}{4}x-5 & \mbox{if } x\in [20,\infty[. \end{cases}
\end{array}$
\end{center}
\noindent The functions $A$, $f$ and $g$ are depicted in Fig. 1. It is obvious that $g$ is convex (but, of course, not strictly convex) and hence super-additive, and it is an easy exercise to show that the (obviously non-linear) function $f$ is sub-additive. From \cite{Greco} it follows that $A^*(x)\le g(x)$ and $A_*(x)\ge f(x)$ for every $x\in [0,\infty[$. Then, clearly, $A^*(x)=g(x)$ for every $x\in [0,4]$, and by the proof of Theorem 1 in {\rm \cite{Sipo}} we have $A^*(x)=g(x)$ also for every $x\in [4,20]$; note that $(20,20)$ is the point of intersection of the lines $y=x$ and $y=\frac{5}{4}x-5$. Moreover, for every $x\in [20,\infty[$ we have $A^*(x)=g(x)=A(x)$ since for every such $x$ the supremum in (\ref{nula}) is equal to $A(x)$, attained for the trivial sum (for $\ell=1$) when $x_1=x$. In an entirely similar way we clearly have $A_*(x)=f(x)$ for every $x\in [0,12]$. To determine the values of $A_*$ for $x\in [12,\infty[$ it is sufficient to realize that every such $x$ can be expressed as a sum of elements lying in the interval $[6,12]$. This together with the fact that $f$ is linear in the interval $[6,\infty[$ and $f(0)=0$ implies that the value of $A_*(x)$ is equal to $f(x)$ for every $x\in [12,\infty[$. Thus, $A_*=f$ and $A^*g$ on $[0,\infty[$.
\eexa

\begin{figure}[h!]\label{fig-non}
  \begin{center}
    \includegraphics[scale=0.45]{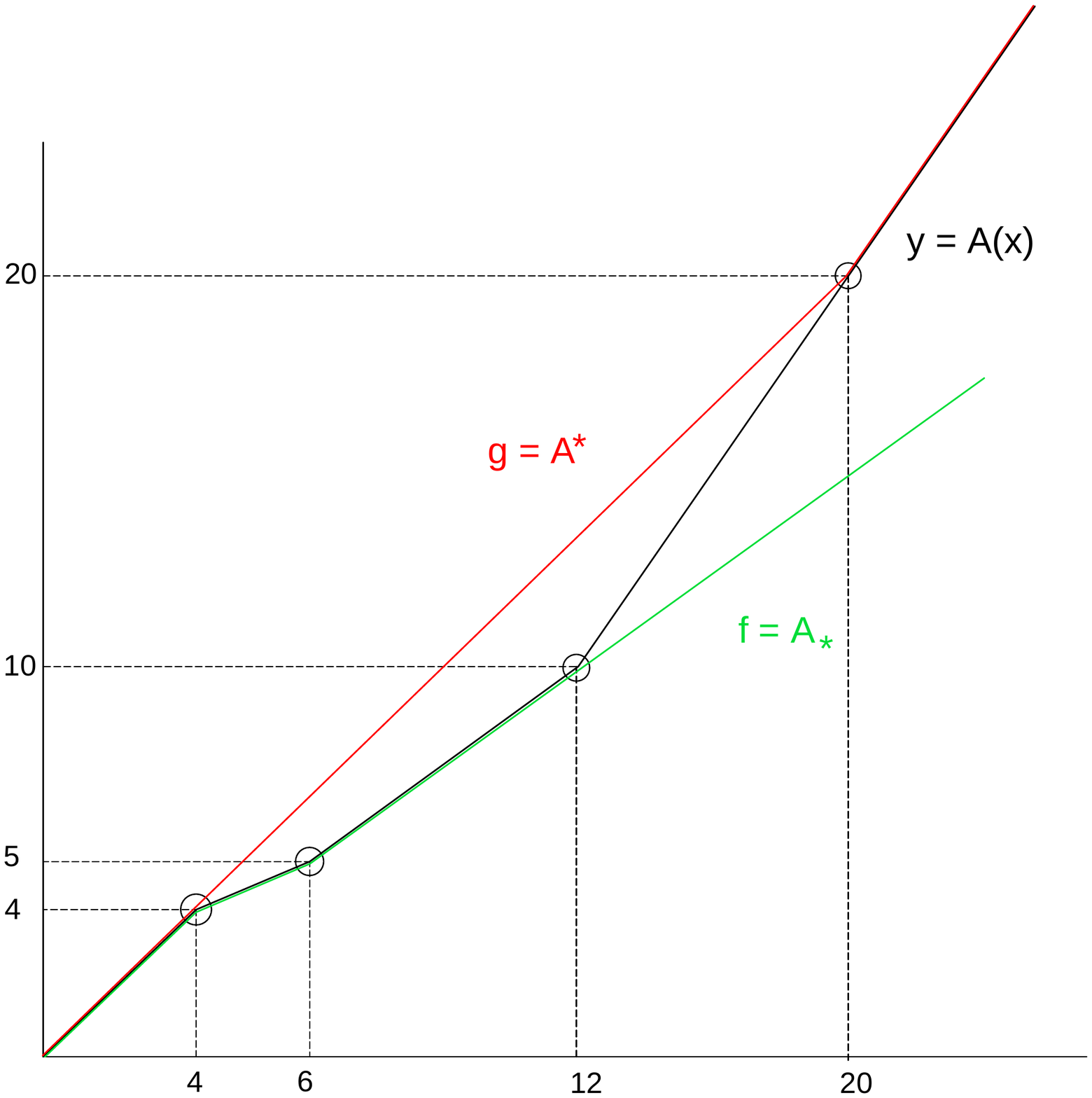}
    \caption{A function $A$ with non-linear $A_*$ and non-strictly convex $A^*$.}
  \end{center}
\end{figure}

This example shows that there is an abundance of pairs of functions $f,g:\ [0,\infty[$ with $f(0)=g(0)$, $f(x)\le g(x)$ for every $x\in [0,\infty[$, such that $f$ is sub-additive and not linear, $g$ is convex (but not strictly convex) and not linear, and yet there exists an aggregation function $A:\ [0,\infty[\to [0,\infty[$ such that $A_*=f$ and $A^*=g$ on $[0,\infty[$. It is clear that one can construct similar examples in the dual version, that is, when $f$ is convex (but not strictly convex) and not linear, and $g$ is super-additive but not linear.
\smallskip

Of course, multi-dimensional examples showing that the assumption of strict directional convexity or concavity cannot be dropped in general can now easily be constructed. A trivial way to do this is to take $\hat f,\hat g:\ [0,\infty[^n\to [0,\infty[$ such that $\hat f(x_1,x_2,\ldots,x_n)=f(x_1)+x_2+\ldots+x_n$ and $\hat g(x_1,x_2,\ldots,x_n)=g(x_1)+x_2+\ldots+x_n$, where $f$ and $g$ are as above. Since $\hat f$ and $\hat g$ have been obtained from $f$ and $g$ by `adding' a linear function in the remaining variables, the preceding analysis carries through and shows that the function $\hat A$ defined by $\hat A(x_1,x_2,\ldots,x_n) = A(x_1)+x_2+\ldots+x_n$ with $A$ as above is an aggregation function on $[0,\infty[^n$ such that $\hat A_*=\hat f$ and $\hat A^*=\hat g$.
\smallskip

It is quite likely, however, that the assumptions on directional convexity and concavity in our main results can still be relaxed in some sensible way. We leave this question as an open problem for further research.

\bigskip

\noindent{\bf Acknowledgement.}~~ The authors would like to than professor R. Mesiar for valuable comments. The first two authors acknowledge support from the APVV 14/0013 and the VEGA 1/0420/15 research grants. Research of the third author was supported by the APVV 0136/12 and the VEGA 1/0007/14 research grants.

\end{document}